\title{Vectors in a Box\thanks{This research
was partially done at the \emph{Gremo Workshop on Open Problems} 2009,
and the support of the ETH Z\"urich is gratefully acknowledged.}}

\newcommand{\cmt}[1]{\ifhmode\newline\fi{\sf *** \ \ #1 \\}}

\date{}

\documentclass[11pt]{article}

\usepackage{a4wide}
\usepackage{latexsym}
\usepackage{amssymb}
\usepackage{amsmath}
\usepackage{graphicx}
\usepackage{epsfig}

\newcommand{\DFG}{Dash, Fukasawa, and G\"unl\"uk}

\author{
{\sc Kevin Buchin}\thanks{Supported by the Netherlands Organisation for Scientific Research (NWO) under project no. 639.022.707}\\
   {\footnotesize Department of Mathematics and Computer Science}\\[-1.5mm]
   {\footnotesize Technical University of Eindhoven, P.O. Box 513}\\[-1.5mm]
   {\footnotesize 5600 MB Eindhoven, Netherlands}\\[-1.5mm]
   {\footnotesize e-mail: {\tt k.a.buchin@tue.nl}}
\and
{\sc Ji\v{r}\'{\i} Matou\v{s}ek}\\
   {\footnotesize Department of Applied Mathematics and}\\[-1.5mm]
   {\footnotesize Institute of Theoretical Computer Science (ITI)}\\[-1.5mm]
   {\footnotesize  Charles University, Malostransk\'{e} n\'{a}m. 25}\\[-1.5mm]
{\footnotesize  118~00~~Praha~1,
   Czech Republic, and}\\
{\footnotesize    Institute of  Theoretical Computer Science}\\[-1.5mm]
{\footnotesize    ETH Zurich,
      8092 Zurich, Switzerland}
\\[-1.5mm] {\footnotesize e-mail: {\tt matousek@kam.mff.cuni.cz}}
\and
{\sc Robin A. Moser}\thanks{This research was partially done
during an internship with Microsoft Research, Redmond, Washington, USA.} \\
{\footnotesize    Institute of  Theoretical Computer Science}\\[-1.5mm]
{\footnotesize    ETH Zurich,
      8092 Zurich, Switzerland}
\\[-1.5mm]   {\footnotesize e-mail: {\tt robin.moser@inf.ethz.ch}}
\and
{\sc D\"om\"ot\"or P\'alv\"olgyi} \\
{\footnotesize    Ecole Polytechnique F\'ed\'erale de Lausanne}
\\[-1.5mm]{\footnotesize    Switzerland}\\[-1.5mm]
{\footnotesize e-mail: {\tt dom@cs.elte.hu}}
}

\newtheorem{theorem}{Theorem}[section]
\newtheorem{definition}[theorem]{Definition}
\newtheorem{prop}[theorem]{Proposition}

\newtheorem{lemma}[theorem]{Lemma}
\newtheorem{corol}[theorem]{Corollary}

\makeatletter
\newcommand{\ProofEndBox}{{\ifhmode\unskip\nobreak\hfil\penalty50 \else
          \leavevmode\fi\quad\vadjust{}\nobreak\hfill$\Box$
            \finalhyphendemerits=0 \par}}
\makeatother
\newcommand{\proofend}{\ProofEndBox\smallskip}

\newcommand\makevec[1]{{\bf #1}}
\def \uu {\makevec{u}}
\def \vv {\makevec{v}}
\def \ww {\makevec{w}}
\let \hastadsa=\aa
\def \aa {\makevec{a}}
\def \bb {\makevec{b}}

\def \xx {\makevec{x}}

\def \zz {\makevec{z}}

\def \rr {\makevec{r}}

\def \ss {\makevec{s}}

\def \cc {\makevec{c}}

\def \zero {\makevec{0}}
\def \one {\makevec{1}}

\newcommand{\R}{{\mathbb{R}}}
\newcommand{\Z}{{\mathbb{Z}}}

\newcommand\eps{\varepsilon}
\newcommand{\conv}{{\mathop {\rm conv} \nolimits}}

\newcommand{\heading}[1]{\vspace{1ex}\par\noindent{\bf #1}}

\long\def\onefigure#1#2{
\begin{figure*}[tbp]
\begin{center}
#1
\end{center}
\caption{#2}
\end{figure*}
}

\newcommand{\labepsfig}[2]  
{\onefigure{\mbox{\includegraphics{#1}}}{\label{f:#1} #2} }

\newcommand{\tgt}{\tau}
\newcommand{\tgti}{\tau_{\pm 1}}

\begin{document}

 \maketitle

\begin{abstract}
For an integer $d\ge 1$, let $\tgt(d)$ be the smallest
integer with the following property:
If $\vv_1,\vv_2,\ldots,\vv_t$ is a sequence of $t\ge 2$
vectors in $[-1,1]^d$ with
$\vv_1+\vv_2+\cdots+\vv_t\in [-1,1]^d$, then
there is a set $S\subseteq \{1,2,\ldots,t\}$
of indices, $2\le|S|\le\tgt(d)$,  such that
$\sum_{i\in S}\vv_i\in [-1,1]^d$.
The quantity $\tgt(d)$ was introduced by \DFG, who showed that
$\tgt(2)=2$, $\tgt(3)=4$, and $\tgt(d)=\Omega(2^d)$,
and asked whether $\tgt(d)$ is finite for all $d$.

Using the Steinitz lemma, in a quantitative
version due to Grinberg and Sevastyanov,
we prove an upper bound of $\tgt(d)\le d^{d+o(d)}$, and based on
a construction of Alon and V\~u, whose main idea goes back
to H\hastadsa stad, we obtain a lower bound of $\tgt(d)\ge d^{d/2-o(d)}$.

These results contribute to understanding
the \emph{master equality polyhedron with multiple rows}
defined by Dash et al., which is a ``universal'' polyhedron
encoding valid cutting planes for integer programs
(this line of research was started
by Gomory in the late 1960s). In particular, the upper bound on
$\tgt(d)$ implies a pseudo-polynomial running time for
an algorithm of Dash et al. for integer programming with a fixed number
of constraints. The algorithm consists in solving a linear program,
and it provides an alternative to a 1981 dynamic programming
algorithm of Papadimitriou.

\end{abstract}

\section{Introduction}

Let $d\ge 1$ and let us consider the unit cube $[-1,1]^d$;
we will call it \emph{the box} in this paper.
We want to construct a large number $t$ of
vectors $\vv_1,\vv_2,\ldots,\vv_t$, each of them lying in the box,
such that
\begin{enumerate}
\item[(i)] the sum $\ss=\vv_1+\vv_2+\cdots+\vv_t$
also lies in the box, but
\item[(ii)] for every \emph{proper} subset
$S\subset [t]$ of indices\footnote{We use the notation
$[t]=\{1,2,\ldots,t\}$.}
with $2\le|S|<t$,
the sum $\sum_{i\in S}\vv_{i}$
lies outside the box (we have to exclude $|S|=1$, since
every $\vv_i$ itself does lie in the box).
\end{enumerate}
So we are interested in long \emph{minimal}
sequences\footnote{Strictly speaking, the order of the vectors
is irrelevant for the considered property, and so one should
perhaps rather speak of sets or multisets of vectors.
However, we find sequences easier to work with for notational
reasons.}
 with sum in the box.  Let $\tgt(d)$ denote the largest $t$
such that a minimal sequence as above exists (it is easy to see
that the definition in the abstract, although phrased
differently, is actually equivalent).

In order to illustrate this definition, let us check that  $\tgt(2)=2$.
We have $\tgt(d)\ge 2$ for all $d$ by definition.
For proving $\tgt(d)\le 2$, we need to show that in every
 sequence $\vv_1,\ldots,\vv_t\in [-1,1]^2$ with sum in the box
there are two vectors with sum in the box.

\labepsfig{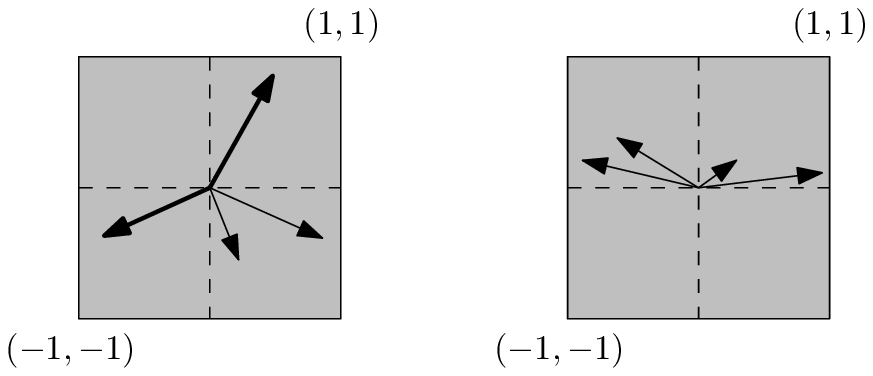}{Illustration to the proof of $\tgt(2)=2$:
Either we find two vectors in opposite quadrants,
or one of the coordinates has the same sign for all the vectors
and can be ignored.}

If two of the vectors lie in opposite quadrants, as in Fig.~\ref{f:tau2}
left, then their sum is in the box and we are done.
Otherwise, some two neighboring quadrants have to be empty,
which means that one of the two coordinates has the same sign
for all the $\vv_i$; w.l.o.g.\ we may assume that
all the $\vv_i$ have a positive $y$-coordinate.
Then the $y$-coordinate can be ignored (since
it lies in $[-1,1]$ for the sum of any subsequence), and it suffices
to show that the sum of some two of the $x$-coordinates
lies in $[-1,1]$. In other words, it now suffices to check
that $\tgt(1)=2$, which we leave to the reader.

An example showing $\tgt(3)\ge 4$ is the sequence
$(1,1,\frac23)$, $(1,-\frac23,-1)$, $(-\frac23,1,-1)$, $(-\frac23,-\frac23,\frac23)$.

The quantity $\tgt(d)$ was introduced by \DFG~\cite{dfg-mepmr-09}
in the context of integer programming (we will discuss the motivation
later). They found the values $\tgt(2)=2$ and $\tgt(3)=4$,
and they asked whether $\tgt(d)$ is finite for all $d$.
We provide a positive answer, with the following upper bound:

\begin{theorem}\label{t:ub}
For all $d\ge 1$ we have
$
\tgt(d)<4(2d)^d.
$
\end{theorem}

Our proof, presented in Section~\ref{s:ub} below,
is based on the so-called Steinitz lemma, in a quantitative version
due to Grinberg and Sevastyanov \cite{GriSev}.

We also show that the upper bound is not far from the truth.

\begin{theorem}\label{t:lb} There is a constant $c>0$ such that
$$
\tgt(d)\ge (cd)^{d/2}
$$
for all $d$ that are powers of~$2$.
\end{theorem}

The proof, given in Section~\ref{s:ub}
 is based on a construction of very ill-conditioned
square matrices with $\pm1$ entries, due to Alon and V\~{u}
\cite{av-ahm-97} (the basic idea going back to
H{\hastadsa}stad \cite{HastadGates}). It seems likely that
the lower bound could be extended to all $d$, instead of
just powers of $2$, but this might need a careful analysis
of another construction from \cite{av-ahm-97}.

There is a natural and, in our opinion, interesting
variant of the quantity $\tgt(d)$, where one again considers
sequences $\vv_1,\ldots,\vv_t\in [-1,1]^d$
satisfying (i) and (ii) above, but the $\vv_i$
are restricted to only $\pm 1$ vectors. Let $\tgti(d)$
denote the corresponding maximum length of such a sequence;
we have $\tgti(d)\le\tgt(d)$ by definition.
We obtain the following slightly weaker lower bound:

\begin{theorem}\label{t:lbi} There is a constant $c>0$ such that
$
\tgti(d)\ge (cd)^{d/4}
$
for all $d$ that are powers of~$2$.
\end{theorem}


\heading{The IP connection. }
The quantity $\tgt(d)$ has been motivated by a connection to
an algorithm for integer programming.

Let us consider an integer program in the form
\begin{equation}\label{e:theIP}
\min\{\cc^T\xx: \xx\in\Z^\ell, A\xx=\bb, \xx\ge\zero\},
\end{equation}
where $A$ is an $m\times \ell$ integer matrix, $\cc\in\Z^\ell$,
and $\bb\in\Z^m$. This optimization problem is well-known
to be NP-hard even for $m=1$, i.e., for a single equality
constraint (this follows, e.g., from the hardness of the knapsack problem).
On the other hand, Papadimitriou \cite{p-cip-81}
proved that if $m$ is fixed and the
entries of $A$ and $\bb$ are small integers, bounded in absolute
value by a parameter $N$, then the integer program can be solved
in \emph{pseudo-polynomial} time. That is, the running
 time can be bounded by a polynomial
in $\ell$ and $N$ (and the input size of $\cc$); the
polynomial depends on~$m$.

Papadimitriou's algorithm is based on dynamic programming
(also see Schrijver~\cite{s-tlip-86} for a description);
it searches for a shortest path in an auxiliary graph.
Dash et al.~\cite{dfg-mepmr-09} provided a completely different
algorithm for the same problem, which
consists in solving a linear program over an auxiliary
polyhedron (the so-called \emph{polaroid}).\footnote{More precisely,
two linear programs are needed. Moreover,
the basic algorithm discussed in \cite{dfg-mepmr-09} solves the
\emph{separation problem} for the set
$Q:=\conv\{\xx\in\Z^\ell, A\xx=\bb, \xx\ge\zero\}$, but a dual version
of it can also be used for the optimization variant.
But here we don't want to go into details.}
They obtained a pseudo-polynomial bound for the number of inequalities
in the linear program, and thus also for the running time,
but only for input integer programs with  $m\le 3$ constraints
(actually, they handled the case of $m=1$ constraint
earlier in \cite{dfg-gmcgp}).
To get pseudo-polynomiality for larger $m$, they needed
the finiteness of $\tgt(m)$. Thus, combined with our
Theorem~\ref{t:ub}, their algorithm  provides an alternative to
Papadimitriou's method.

We won't review the algorithm here; we only recall some of the
key concepts, and then we indicate how $\tgt(m)$ is related
to the linear program.

The approach of Dash et al. goes back to a paper of Gomory
\cite{g-sprcp-69}. In that remarkable work, which introduced
several important ideas of modern polyhedral combinatorics,
Gomory defined a certain ``universal'' polyhedron,
the \emph{master cyclic group polyhedron}, whose faces encode
all instances of integer programs in a certain class
(see, e.g., \cite{aardal-al} for an introduction).
Dash et al.~\cite{dfg-mepmr-09} use the somewhat related
concept of the  \emph{master equality polyhedron} $K^m(N,\bb)$,
which we recall below. For $m=1$, it was introduced
in an earlier paper by Dash et al.~\cite{dfg-gmcgp},
who attribute its origin to a 2005 talk of Uocha,
and it contains as a face Gomory's master cyclic group polyhedron,
as well as the \emph{master knapsack polyhedron} of Ara\'oz.

Let
$$I:=\{-N,-N+1,\ldots,N\}^m,
$$
and let $\bb\in I$ be a vector corresponding to the right-hand side in
(\ref{e:theIP}). Then  the master equality polyhedron
resides in $\R^I$ and it is defined as
$$
K^m(N,\bb):=\conv\{\xx\in\Z^I: \sum_{\vv\in I} x_\vv\vv=\bb,\xx\ge\zero\}.
$$
It turns out that the separation problem
for the integer program (\ref{e:theIP}) can be reduced to
the separation problem for $K^m(N,\bb)$.

For $m=1$, Dash et al.~\cite{dfg-gmcgp} obtained a description of
a \emph{nontrivial polar} $T$ of $K^1(N,b)$, i.e. a polyhedron
whose vertices correspond to the nontrivial facets of $K^1(N,b)$,
and thus reduced the separation problem for $K^1(N,b)$
to optimizing a linear function over $T$. It is important that $T$
is described by polynomially many linear constraints.

In \cite{dfg-mepmr-09} Dash et al.\ describe an example showing
that for $m\ge 2$, a sufficiently compact description
of a nontrivial polar of $K^m(N,\bb)$ may be very hard to find,
if it exists at all.  However, they defined
the new notion of a \emph{polaroid} of $K^m(N,\bb)$,
and proved that the separation problem
for $K^m(N,\bb)$ can be reduced to a suitable linear
program over the polaroid. Again, one needs a polynomial bound on
the number of constraints in this linear program,
and here $\tgt(m)$ enters the game.

The variables of the linear program are $\pi_{\vv}$ for all integer vectors
$\vv\in I$, and the main kind of constraints
in it are \emph{subadditivity constraints} of the form
\begin{equation}\label{e:subad}
\pi_{\vv_1}+\cdots+\pi_{\vv_t}\ge \pi_{\vv_1+\cdots+\vv_t},
\end{equation}
where $\vv_1,\ldots,\vv_t\in I$ are vectors whose sum also
lies in $I$. Now suppose, for example,
that $\vv_1+\vv_2+\cdots+\vv_k$ also lies in $I$, for some
$k$ with $2\le k<t$. Then (\ref{e:subad}) is a consequence
of the subadditivity constraints
$$
\pi_{\vv_1}+\cdots+\pi_{\vv_k}\ge \pi_{\vv_1+\cdots+\vv_k}
\ \ \mbox{and}\ \
\pi_{\vv_1+\cdots+\vv_k}+\pi_{\vv_{k+1}}+\cdots+\pi_{\vv_t}\ge
\pi_{\vv_1+\cdots+\vv_t}.
$$
Similarly, whenever there is an $S\subset[t]$ with $2\le|S|<t$ and
$\sum_{i\in S}\vv_i\in I$, the constraint (\ref{e:subad}) is implied
by subadditivity constraints with a smaller number of terms.
Thus, it is sufficient to consider only subadditivity
constraints with $t\le\tgt(m)$.

The quantity $\tgt(m)$ gives only an upper bound on the number
of non-redundant subadditivity constraints.
Dash et al.~\cite{dfg-mepmr-09} define another quantity
$k^*(m)\le \tgt(m)$,
which is related to the number of non-redundant constraints
more directly.
In Corollary~\ref{c:k*} we will give a lower bound of
$(cm)^{m/4}$ for $k^*(m)$,
which shows that it is not much smaller than $\tgt(m)$.

\section{The Upper Bound}\label{s:ub}

Let $B\subset\R^d$ be a $d$-dimensional closed convex body symmetric
about $\zero$ (in other words, the unit ball of a norm on $\R^d$).
Riemann and L\'evy in the 19th century raised the question of whether
there exists a number $m=m(B)$, depending only on $B$, such that
the vectors of an arbitrary finite set (or multiset) $V\subset B$
with $\sum_{\vv\in V}\vv=\zero$
can be ordered into a sequence $\vv_1,\vv_2,\ldots,\vv_n$
so that each of the partial sums $\vv_1+\vv_2+\cdots+\vv_k$,
$1\le k\le n$, belongs to the expanded body $mB=\{m\xx:\xx\in B\}$.
The first complete proof of a positive answer was given by
Steinitz \cite{SteinitzLemma}. The strongest known
quantitative version, with $m(B)=d$ for all bodies in $\R^d$,
is due to Grinberg and Sevastyanov \cite{GriSev}
(their beautiful proof can also be found in B\'ar\'any's survey
\cite[Theorem 2.1]{b-pld-08}).

\begin{theorem}[Steinitz lemma]
\label{t:Steinitz}
Let $B$ be a symmetric convex body in $\R^d$,
and let $V \subset B$ be a finite set (or multiset) of vectors
satisfying $\sum_{\vv \in V}\vv=0$.
Then there is an ordering
$\vv_1,\vv_2, \mathellipsis, \vv_n$ of the elements of $V$ such that for all
$k=1,2,\ldots,n$, we have $\sum_{i=1}^k \vv_i \in d B$.
\end{theorem}

\heading{Proof of Theorem~\ref{t:ub}. }
Given a collection of $t$ vectors in the box whose sum also lies in the box,
and assuming $t\ge t_0:=4(2d)^d$, we want to find a proper subcollection
with sum in the box.

Here it will be more convenient to regard the given collection
of $t$ vectors as a multiset $W$ (rather than as a sequence)---later
we will obtain suitable ordering of the vectors in $W$ from the
Steinitz lemma.

Thus, $W$ is a multiset of $t$ vectors; we let
$\ss:=\sum_{\vv\in W}\vv\in [-1,1]^d$ be their sum.
We apply
the Steinitz lemma as above with $B=[-1,1]^d$
and $V:=W\cup\{-\ss\}$ (this is again a multiset,
with $t+1$ vectors). This yields an ordering
$$
\vv_1,\vv_2,\ldots,\vv_i,-\ss,\vv_{i+1},\ldots,\vv_t
$$
such that the sum of the first $k$ terms lies in $dB$ for every
$k=1,2,\ldots,t+1$.

First, let us assume that the ``artificial'' element $-\ss$ falls
in the second half of the above sequence; that is,
$i\ge t/2 \ge t_0/2$. Let us subdivide
the blown-up cube $dB$, whose side length is $2d$, into
$(2d)^d$ cubes of side $1$, of the form $[0,1]^d+\zz$,
$\zz\in\{-d,-d+1,\ldots,d-1\}^d$.

Let $\ss_k:=\vv_1+\cdots+\vv_{k}$ be the $k$th partial sum,
$k=0,1,\ldots,t_0/2$, and let us consider
every second of these, i.e., the points
$\ss_0,\ss_2,\ss_4,\ldots,\ss_{t_0/2}$.
These are more than $t_0/4=(2d)^d$ points in $dB$,
and so some two of them, $\ss_{2i}$ and $\ss_{2j}$,
$0\le i<j\le t_0/4$, fall in the same cube $[0,1]^d+\zz$.

Then the difference $\ss_{2j}-\ss_{2i}$ lies in the box.
At the same time, it equals
$$
\vv_{2i+1}+\vv_{2i+2}+\cdots+\vv_{2j},
$$
and so we have found the desired proper submultiset of $W$
with sum in the box (and with at least two elements).

It remains to deal with the case when the artificial
element $-\ss$ lies in the first half of the sequence---then
we use the same kind of argument as above for the second half.
Theorem~\ref{t:ub} is proved.
\proofend



\section{Lower bounds}\label{s:lb}


All of our lower bounds are based on results
of Alon and V\~{u}  \cite{av-ahm-97}.
The main theme of that paper are ill-conditioned
matrices with $\pm 1$ entries (or $0/1$ entries; these
two settings are not very different).

Let us consider a nonsingular $d\times d$ matrix $A$
whose entries are $+1$'s and $-1$'s, and let $\chi(A)$
be the maximum of the absolute values of the entries
of $A^{-1}$; the larger $\chi(A)$, the more ill-conditioned the
matrix~$A$ is. Let $\chi(d):=\max_A\chi(A)$, where
the maximum is over all $d\times d$ nonsingular $\pm 1$ matrices.

Alon and V\~{u} showed that $\chi(d)=d^{d/2+o(d)}$
and gave several interesting applications.
The main achievement was the (surprisingly large)
lower bound, which was obtained by an explicit construction.
For that purpose, Alon and V\~{u} modified and extended
a construction of H{\hastadsa}stad \cite{HastadGates},
which was formulated in a different setting,
namely, in the language of threshold gates.

In this section we will provide three lower bound constructions
for vectors in the box. The first construction is the simplest:
It follows rather directly from a result explicitly stated
in \cite{av-ahm-97}, but it loses a factor of two in the dimension,
leading only to the lower bound of $\tgt(d)\ge d^{d/4-o(d)}$.
Then we give another, different and more complicated construction,
which gives $\tgt(d)\ge d^{d/2-o(d)}$ and proves Theorem~\ref{t:lb},
and finally, we modify the latter construction so that
only $\pm1$ vectors are used, obtaining Theorem~\ref{t:lbi}.

\heading{Large quantities. } In order to avoid boring formulas
or repetitive phrases, we introduce the following piece
of terminology. If $Q=Q(d)$ is some quantity in the forthcoming
constructions, depending on $d$, we say that $Q$ is
\emph{large} if there is a constant $c>0$
such that $|Q|\ge (cd)^{d/2}$ holds for all of the
considered values of $d$ (in the first construction
we will consider all $d$, while in the second and third
only powers of two). In particular, we note that a large
quantity can be either positive or negative.

We note that being large is immune to division by exponential
factors; e.g., if some $Q$ is large, then $Q/2^d$
is large as well.


\subsection{The First Construction}

Here we use the following result of Alon and V\~u:

\begin{prop}[{\cite[Proposition~3.4.3]{av-ahm-97}}]
\label{pro:alonvufirst}
For every $d\ge 1$, there exists a $d \times (d+1)$
matrix $C$ of rank $d$ with entries $+1$ and $-1$
such that any nonzero integral solution $\zz$ of the system
$C\xx=\zero$
has at least one large component.
\end{prop}

\heading{The construction. }
The following construction will show that $\tgt(2d)$ is large (for all $d$),
or in other words, that $\tgt(d)\ge (cd)^{d/4}$.

Let $C$ be the $d\times(d+1)$ matrix as in
Proposition~\ref{pro:alonvufirst}. Since the system
$C\xx=\zero$ is homogeneous, with rational coefficients,
 and has fewer equations than
unknowns, there exist nonzero integral solutions.

Let $\zz\in\Z^{d+1}$ be a nonzero integral solution of $C\xx=\zero$
with the smallest possible $L_1$ norm, i.e.,
minimizing $\|\zz\|_1=\sum_{j=1}^{d+1}|z_j|$.
Let us set $t:=\|\zz\|_1$; by the above, $t=t(d)$ is large.

After possibly flipping the signs of some of the columns of $C$,
we may assume that all components of $\zz$ are nonnegative.
Let $\cc_j$ denote the $j$th column of $C$, and let
$\ww_1,\ww_2,\ldots,\ww_t$ be an (auxiliary) sequence of vectors containing
$z_j$ copies of each $\cc_j$, $j=1,2,\ldots,d+1$.

The $\ww_i$ are vectors in the box (even $\pm1$ vectors),
and we have $\sum_{i=1}^t\ww_i=\sum_{j=1}^{d+1}z_j\cc_j=C\zz=\zero$.

It is also easy to see that no proper subsequence of the $\ww_i$
has sum in the \emph{interior} of the box. Indeed, the sum of
any subsequence is an integral vector, so if it lies in
$(-1,1)^d$, it has to be $\zero$. But choosing a proper subsequence
of the $\ww_i$ corresponds to choosing multiplicities
$z'_1,\ldots,z'_{d+1}$, with $z'_j\le z_j$ for all $j$
and with at least one of the inequalities strict.
So a proper subsequence with zero sum corresponds
to a nontrivial solution $\zz'$ of $C\zz'=\zero$
with $\|\zz'\|_1<t=\|\zz\|_1$, contradicting the assumed minimality
of $\|\zz\|_1$.

Thus, the $\ww_i$ almost achieve what we want, but only almost,
since there may be some sums of proper subsequences on the boundary
of the box. We get around this by a simple dimension-doubling trick.
\medskip

Let us set $\eps := 1/(10t)$, say. Let $\ww'_i$ be the vector
obtained from $\ww_i$ by replacing all $-1$ components
by $-(1-\eps)$ and keeping all $+1$ components. Similarly,
$\ww''_i$ is obtained by keeping the $-1$ components of $\ww_i$
and replacing $+1$'s by $1-\eps$. Thus, for example,
if we had $\ww_i=(+1,+1,-1)$,  then $\ww'_i=(+1,+1,-(1-\eps))$
and $\ww''_i=(1-\eps,1-\eps,-1)$. Finally,
let $\vv_i\in\R^{2d}$ be obtained by concatenating $\ww'_i$
and $\ww''_i$.

We claim that this sequence $\vv_1,\ldots,\vv_t$ witnesses
$\tgt(2d)\ge t$. Clearly, the $\vv_i$ lie in the box.
Moreover, since all the $\ww_i$ sum to $\zero$ and
$\|\ww_i-\ww'_i\|_\infty\le\eps$, $\|\ww_i-\ww''_i\|_\infty\le\eps$,
we have $\sum_{i=1}^t \vv_i\in [-t\eps,t\eps]^{2d}\subset [-1,1]^{2d}$.

Next, let us consider a proper subset $S\subset[t]$,
$2\le|S|<t$. We already know that $\sum_{i\in S}\ww_i\ne\zero$;
let us fix a coordinate $k$ in which this sum has a nonzero
component. Let $a$ be the number of $+1$'s in the $k$th coordinate
of the sum, and let $b$ be the number
of $-1$'s there; that is,
$a:=|\{i: (\ww_i)_k=1\}|$, $b:=|S|-a=|\{i: (\ww_i)_k=-1\}|$.
We have $a\ne b$.

Then $\left(\sum_{i\in S}\vv_i\right)_k=
\left(\sum_{i\in S}\ww'_i\right)_k=(a-b)+b\eps$
and $\left(\sum_{i\in S}\vv_i\right)_{d+k}=
\left(\sum_{i\in S}\ww''_i\right)_k=(a-b)-a\eps$; we claim
that at least one of these numbers falls outside $[-1,1]$.
Indeed, this is clear if $|a-b|\ge 2$. If $a-b=1$,
then $b\ge 1$ (since $a+b=|S|\ge 2$),
and so $(a-b)+b\eps\ge 1+\eps>1$. Similarly, for $a-b=-1$ we find that
$a-b-a\eps\le -1-\eps<-1$. Thus, $\sum_{i\in S}\vv_i$ is not in the
box and Theorem~\ref{t:lb} is proved.
\proofend

\heading{A lower bound for the quantity \boldmath$k^*(m)$. }
As was mentioned in the last part of the introduction,
Dash et al.~\cite{dfg-mepmr-09} define an integer function
$k^*(m)$, which is bounded above by $\tgt(m)$,\footnote{They don't
prove this explicitly, but it can be seen from the proofs
of their Theorems~4.4 and 4.7.} but which is more
directly related to the number of constraints in their linear program.
Here we won't recall the definition of $k^*(m)$, since we won't
use it directly. Rather, we will rely  on a
property of $k^*(m)$, which is expressed in Lemma~6.1 of \cite{dfg-mepmr-09},
and which in our notation can be re-phrased as follows.

\begin{lemma}[Dash et al. \cite{dfg-mepmr-09}]\label{l:k*}
 Let $m\ge 1$, and suppose that there
are rational  vectors $\vv_1,\ldots,\vv_t\in [-1,1]^{m}$,
with $\ss:=\sum_{i=1}^t\vv_i\in[-1,1]^d$, such that $\sum_{i\in S}\vv_i$ lies
outside the box for every $S\subset[t]$, $2\le|S|<t$,
and moreover, for every choice of nonnegative integer coefficients
$q_1,\ldots,q_t$ with $1\le q:=\sum_{i=1}^t q_i\le t/2$,
the vector $\ss-\ss'$, where $\ss':=\sum_{i=1}^tq_i\vv_i$,
 also lies outside the box. Then $k^*(m)\ge t$.
\end{lemma}

It turns out that the vectors $\vv_1,\ldots,\vv_t$ constructed in the previous
proof also have the  additional property in the above lemma,
and so we get $k^*(m)\ge (cm)^{m/4}$ for all even $m$:

\begin{corol}\label{c:k*} We have $k^*(m)\ge (cm)^{m/4}$ for all even
$m$, where $c$ is a positive constant.
\end{corol}

\heading{Sketch of proof. } Let us set $m=2d$ and use the vectors
$\vv_1,\ldots,\vv_t\in[-1,1]^{2d}$ from the previous proof. The only property
which we haven't yet verified for them is
 the ``moreover'' part in Lemma~\ref{l:k*}, and we do this now.
We need to assume $t\ge 4$ (which we can since the bound is asymptotic
and so very small values of $m$ can be ignored).

It is easily checked that the first $d$ coordinates of $\ss$
equal $\frac12t\eps$  and the remaining $d$ coordinates are
$-\frac12t\eps$.

Arguing as in the previous proof,
we get that  $\sum_{i=1}^tq_i\ww_i\ne\zero$ (the fact that one
$\ww_i$ may appear several times in the sum makes no difference),
and so there is some coordinate $k$ where the number $a$
of $+1$ contributions differs from the number $b$ of $-1$ contributions,
$a+b=q$.
(More formally, $a=\sum_{i:(\ww_i)_k=+1}q_i$,
$b=\sum_{i:(\ww_i)_k=-1}q_i$.) Let us suppose, for example, that
$a<b$; then we calculate that the $k$th coordinate of
$\ss-\ss'$ equals $\frac12t\eps-(a-b)-b\eps$, which is surely
above $+1$ for $a-b\le -2$. For $a-b=-1$ it equals
$1+(t/2-b)\eps$, and we have $b=\frac {q+1}2\le
\frac{t/2+1}2<t/2$. For $a>b$ we argue similarly using
the $(k+d)$th coordinate.
\proofend

\subsection{The Second Construction}

The stronger construction used in the proof of Theorem~\ref{t:lb}
is based on the following result of Alon and V\~{u}
\cite{av-ahm-97}.

\begin{prop}
\label{pro:alonvusecond}
For every $d$ that is a power of $2$, there exists
a $d \times d$ nonsingular matrix $A$ with entries $\pm 1$ such that
the matrix $B:= 2^d A^{-1}$ is integral, has nonnegative row sums,
and all the entries in the first row of $B$ are nonnegative and large.
\end{prop}

This statement is not explicitly formulated in \cite{av-ahm-97};
rather, it can be combined from several remarks
scattered throughout that paper, so we recall a
(very easy) proof from more a explicit statement
in~\cite{av-ahm-97}.

\heading{Proof.}
Let $d$ be a power of $2$. In the proof of Theorem 2.1.1 in
\cite{av-ahm-97}, Alon and V\~{u} construct an $d \times d$
nonsingular matrix $\tilde A$ with entries $\pm 1$ such that
there exists a column of $\tilde A^{-1}$
in which all entries are large.

They also show that $\det(\tilde A)=2^{d-1}$. By transposing $\tilde A$
and reordering its columns, we can guarantee
that the first row of the inverse matrix consists
of large entries. Since changing the sign of a column
changes the sign of the corresponding row of the inverse matrix,
by flipping the signs of suitable rows we can make
all entries in the first row of the inverse nonnegative.
Finally, by flipping the signs of some columns
we can arrange for nonnegativity of the row sums of the inverse.
In this way we obtain the desired~$A$.

 Since all the operations performed above
preserve the determinant, we still have $\det(A)=2^{d-1}$, and since
the adjoint ${\rm adj}(A)$ is integral, $B=2^{d}A^{-1}$
is integral as well.
\proofend

\medskip

For a vector $\xx\in\R^n$, the notation $\xx>\zero$ means that all
entries of $\xx$ are nonnegative and $\xx\ne\zero$.

\begin{corol}
\label{cor:alonvuproperty}
Let $A$ be as in the previous proposition and
let $\bb>\zero$ be an integral vector.
Then the (unique) solution $\zz$ of $A\xx=\bb$ has the first
component $z_1$ positive and large.
\end{corol}

\heading{Proof.}
Since $\zz=A^{-1}\bb$, $z_1$ is a linear combination of the
entries of the first row of $A^{-1}$ with nonnegative
integer coefficients (given by $\bb$), at least one of them nonzero.
Since the entries in the first row of $A^{-1}$ are all large
and positive, the corollary follows.
\proofend

\medskip




\heading{Proof of Theorem~\ref{t:lb}. }
Using the matrix  $A$ provided by
Proposition~\ref{pro:alonvusecond}, we construct a sequence
$\vv_1, \mathellipsis, \vv_t$ of vectors in $[-1,1]^{d+1}$ with
$t=t(d)$ large and with sum in the box. This time we won't show
that the sequence is minimal; rather, we will prove that
every subsequence with at least $2$ terms and with sum in the box
has to have a large number of terms.

Let $r_j$
be the sum of all entries in the $j$th row of $B=2^dA^{-1}$, which, as
Proposition~\ref{pro:alonvusecond} asserts, is a nonnegative
integer. We can also write $\rr=(r_1,\ldots,r_d)=B\one$,
where $\one$ is the all 1's vector.
Let $R:=\sum_{j=1}^d r_j$, and let $\alpha:=\frac{2^d}R$.
We note that since $R$ is large, may assume $\alpha\in (0,1)$.

Let $\aa_j$ denote the $j$th column of $A$, and let
$\overline{\aa}_j=(\aa_j,1)\in\R^{d+1}$ be
obtained from $\aa_j$ by appending the component $1$ to the end.

Our sequence $\vv_1, \mathellipsis, \vv_t$ consists of $r_i$ copies
of $\overline{\aa}_j$, $j=1,2,\ldots,d$, and of $R$ copies of the
vector $\cc:=(-\alpha\one,-1)\in\R^{d+1}$, making a
total number of $t = 2R$ vectors.

The vector $\cc$ and its multiplicity were chosen so that
the sum of all vectors in the sequence is $\zero$, as we now check.
For the $(d+1)$st coordinate this is equivalent to
$\sum_{j=1}^n r_j=R$. The vector consisting of the first
$d$ entries of $\sum_{i=1}^t\vv_i$ equals
$$
\biggl(\sum_{j=1}^n r_j\aa_j\biggr) - R\alpha\one=
A\rr-2^d\one=
A(B\one) -2^d\one= 2^d AA^{-1}\one-2^d\one=\zero.
$$

It remains to show that if for some  $S\subseteq [t]$, $|S|\ge 2$,
the sum $\sum_{i\in S}\vv_i$ lies in the box, then
$|S|$ is large.

Choosing a subsequence corresponds to choosing multiplicities
of the vectors $\overline{\aa}_1,\ldots,\overline{\aa}_d$ and $\cc$;
we denote these
multiplicities by $z_1,\ldots,z_d$ and $k$, respectively.
The number of terms is $\sum_{j=1}^d z_j + k \ge 2$.

Let us suppose that the sum $\ss':=\sum_{i\in S}\vv_i$ lies
in $[-1,1]^{d+1}$.
First we check the $z_j$ can't be all $0$. If we had $\zz=\zero$,
then we would get $k\ge 2$, and
$s'_{d+1}\le -2$---a contradiction. So $\zz>\zero$.

Similarly we find, using the last coordinate again, that $k>0$.
Indeed, if we had $k=0$, then $s'_{d+1}=\sum_{j=1}^d z_j\ge 2$.
Thus, $k>0$ as claimed.

The vector of the first $d$ coordinates of $\ss'$
equals $A\zz-k\alpha\one$.
We consider the vector $\bb:=A\zz$. Clearly, it is integral
and nonzero (since the only solution of $A\xx=\zero$ is $\zero$,
while $\zz\ne\zero$). We claim that $\bb\ge\zero$.
Indeed, if $b_j<0$, then $b_j\le -1$ by integrality,
and so we would get
$s'_j=b_j-k\alpha\le -1-\alpha<-1$ (using $k\ge 1$)---a
contradiction.

We have shown that $A\zz=\bb$ with $\bb>\zero$ integral, and
we can apply Corollary~\ref{cor:alonvuproperty} to conclude that
$z_1$ is large. This also means that $|S|$ is large.
\proofend


\subsection{The Third Construction: \boldmath $\pm1$ Vectors}

Here we will prove Theorem~\ref{t:lbi}, the lower bound
for $\tgti(d)$. To this end, we will exhibit a
a sequence $\vv_1,\ldots,\vv_t$, $t$ large, of $\pm1$ vectors in $\R^{2d+1}$,
with sum $\zero$ and such that
every subsequence of length at least $2$
with sum in the box has a large number of terms.

As in the previous subsection, we use the matrix $A$ provided by
Proposition~\ref{pro:alonvusecond}, we set $r_j := (2^d A^{-1} \one)_j$,
$R := \sum_{j=1}^n r_j$, and $t = 2R$.
By Corollary~\ref{cor:alonvuproperty}, $t$ is large. Moreover, we note
that $R$ is divisible by two because, as we demonstrated in the
previous section, $\sum_{j=1}^d r_j \aa_j = 2^{d} \one$, and since all
the $\aa_j$ are vectors with $\pm 1$ entries, we need an even number of
them to reach a point with even coordinates. Therefore, $t$ is
divisible by~$4$.

Now let $\uu_1, \mathellipsis, \uu_t \in \{+1,-1\}^d$ be a sequence
of vectors that consists of $r_i$ copies of the $i$th column of $A$.
We build the vectors  $\vv_1, \mathellipsis, \vv_t \in \{+1,-1\}^{2d+1}$
as follows. For $i=1,2,\ldots,t/2$ and $j\in[2d+1]$, we let
$$ (\vv_i)_j \; := \;
   \begin{cases}
     (\uu_i)_j & \text{ if $1 \le j \le d$} \cr
     +1       & \text{ if $d+1 \le j \le 2d$ and $1 \le i \le t/4-2^{d-1}$} \cr
     -1       & \text{ if $d+1 \le j \le 2d$ and $t/4-2^{d-1} < i \le t/2$} \cr
     +1       & \text{ if $j=2d+1$.}
   \end{cases}
$$
Then for $i=t/2+1,\ldots,t$ and $j\in[2d+1]$, we let
$$ (\vv_i)_j \; := \;
   \begin{cases}
     (\vv_{i-t/2})_{j+d} & \text{ if $1 \le j \le d$} \cr
     (\vv_{i-t/2})_{j-d} & \text{ if $d+1 \le j \le 2d$} \cr
     -1       & \text{ if $j=2d+1$.}
   \end{cases}
$$

We first claim that the sequence sums to $\zero$.
Just as in the last section,
$\sum_{i=1}^{t/2} \uu_i = 2^{d} \one$.
Therefore, by the definition of the $\vv_i$,
$$ \sum_{i=1}^{t/2} \vv_i \; = \;
   (\underbrace{2^d,2^d, \mathellipsis, 2^d}_{d \text{ times}},
    \underbrace{-2^d, -2^d, \mathellipsis, -2^d}_{d  \text{ times}},
    t/2
   ), $$
and
$$ \sum_{i=t/2+1}^{t} \vv_i \; = \;
   (\underbrace{-2^d,-2^d, \mathellipsis, -2^d}_{d \text{ times}},
    \underbrace{2^d, 2^d, \mathellipsis, 2^d}_{d  \text{ times}},
    -t/2
   ). $$
In conclusion, the total sum is zero.

Now let us consider an index set $S \subseteq [t]$, $|S|\ge 2$,
and let us suppose that
 $\ss := \sum_{i \in S} \vv_i \in \{-1,0,1\}^d$, and that
$|S|$ is not large. By analyzing several cases,
we will show that this leads to a contradiction.

Let $S_1 := S \cap [t/2]$ and $S_2 := S \setminus S_1$.
 Let $\ss_1 := \sum_{i \in S_1} \vv_i$
and $\ss_2 := \sum_{i \in S_2} \vv_i=\ss-\ss_1$.
Moreover, let $\ss_1^{(1)} \in \{-1,1\}^d$ and
$\ss_1^{(2)} \in \{-1,1\}^d$ be the projections of $\ss_1$ onto
the coordinates $1$ through $d$
and the coordinates $d+1$ through $2d$, respectively,
and similarly for $\ss_2^{(1)},\ss_2^{(2)}$.

First let us
suppose  that $(\ss_1)_{2d+1}=0$. Then $|S_1|=0$,
 and since $|S|\ge2$, we have $|S_2|\ge2$. But
then $(\ss)_{2d+1} = (+1)|S_1| + (-1)|S_2| < -1$---a contradiction.
Therefore, $(\ss_1)_{2d+1}>0$.
Symmetrically, $(\ss_2)_{2d+1} < 0$.

Next, let us suppose
that $\ss_1^{(1)}=\zero$. Since $|S_1|>0$ and
 $\ss_1^{(1)}$ is a linear
combination of the columns of $A$, we would get that
$A \xx = \zero$ has a nonzero integral solution, which is not
the case, and so we can conclude $\ss_1^{(1)} \ne \zero$.
Symmetrically, $\ss_2^{(2)} \ne \zero$.

Now we suppose that $\ss_1^{(2)} \ne \zero$.
By the way that vector is composed, we then have
$\ss_1^{(2)} = k \one$ for some nonzero integer $k$.
First let us assume $k<1$. Then, since $\ss$ lies
in the box,   $\ss_2^{(2)}\ge\zero$, and
since we have shown $\ss_2^{(2)} \ne \zero$,
we even have $\ss_2^{(2)} > \zero$. According to
Corollary~\ref{cor:alonvuproperty}, that implies that
 $\ss_2^{(2)}$ is a sum of a large number copies of the columns
of $A$, and in this case $|S|$ is large, contrary to
the assumption.

So we may suppose $k>1$.
For analogous reasons, this implies that $\ss_2^{(2)}<\zero$.
But then, again by Corollary~\ref{cor:alonvuproperty},
it is impossible to express $\ss_2^{(2)}<\zero$ as a nonnegative
integer linear combination of the columns of $A$
(since the corollary implies that for $\bb<0$,
an integral solution of $A\xx=\bb$ has a negative component)---a
contradiction.

Having dealt with the case $\ss_1^{(2)}\ne \zero$,
we now assume $\ss_1^{(2)}= \zero$; symmetrically, we may
assume $\ss_2^{(1)} = \zero$ as well.

Now the total sum $\ss$ has
$\ss_1^{(1)} \ne \zero$ in the first $d$ coordinates
and $\ss_2^{(2)} \ne \zero$ in the coordinates
$d+1$ through  $2d$.
 For parity reasons we get
 $\ss_1^{(1)}, \ss_2^{(2)} \in \{-1,1\}^d$. This
implies that $|S_1| \equiv |S_2| \equiv 1\,({\rm mod}\,2)$.
 On the other hand, we have shown
$\ss_1^{(2)}=\ss_2^{(1)} = \zero$,
 and this gives $|S_1| \equiv |S_2| \equiv 0\,({\rm mod}\,2)$---a
contradiction.
Theorem~\ref{t:lbi} is proved.
\proofend

\section{Conclusion}

It would be interesting to determine the asymptotics of $\tgt(d)$
more precisely. We conjecture that the truth should be close
to the lower bound, i.e., of order roughly $d^{d/2}$.

One way of improving on the upper bound might be to get a factor
better than $d$ in the Steinitz lemma (Theorem~\ref{t:Steinitz}) for the case
$B=[-1,1]^d$. It is known, and
not hard to see, that if $B$ is the unit ball of the
$\ell_1$ norm, then the factor cannot be smaller than $\frac d2$.
However,
it is possible that the factor of $O(\sqrt d)$ suffices
for $B=[-1,1]^d$ (from which
the bound $d^{d/2+o(d)}$
for $\tgt(d)$ would follow).

However, apparently there is no improvement over $d$ for any $B$ known,
and the problem may be hard. As B\'arany
\cite{b-pld-08} puts it, for the case where $B$ is the
Euclidean ball, ``even the much weaker $o(d)$ estimate
seems to be out of reach though quite a few mathematicians
have tried,'' and for $B=[-1,1]^d$ ``there is no proof in sight
even for the weaker $o(d)$ estimate.''

\subsection*{Acknowledgment}

We would like to thank Tibor Szab\'o for raising the problem
at the GWOP'09 workshop, Sanjeeb Dash for prompt answers to our
questions, and Patrick Traxler for useful discussions.

\bibliographystyle{plain}
\bibliography{vecbox}

\end{document}